\documentclass[12pt]{article}

\usepackage{amssymb,latexsym   }
\usepackage[centertags]{amsmath}

\newtheorem{thm}[equation]{Theorem}
\newtheorem{lemma}[equation]{Lemma}

\newtheorem{cor}[equation]{Corollary}

\newtheorem{claim}[equation]{Claim}

\numberwithin{equation}{section}

\newcommand{\C}{\mathbb{C}}
\newcommand{\D}{\mathbb{D}}
\newcommand{\clD}{\overline{\mathbb{D}}}
\newcommand{\T}{\mathbb{T}}
\newcommand{\A}{\operatorname{Area}}
\newcommand{\Zf}{\mathcal{Z}_f}

\title{Coarse equidistribution of the argument of entire functions of
finite order}

\author{F\"edor Nazarov and Mikhail Sodin\thanks{Supported by
the Israel Science Foundation of the Israel Academy of Sciences
and Humanities.}}

\begin{document}
\maketitle

\centerline{\em To Paul Koosis and Iossif Ostrovskii on occasion of
their birthdays}

\bigskip

\section{Main results}

There is an impressive amount of classical results on the
asymptotic  behaviour of the absolute value of entire functions
(Phragm\'en-Lindel\"of-type theorems, minimum modulus theorems,
and their numerous ramifications).
All of them are based on subharmonicity of the function $\log
|f|$. At the same time, very little is known about the argument of
entire functions. Here we present several results (motivated by
\cite[Theorem~2.2]{NPS}) that show somewhat surprising
equidistribution patterns in the asymptotic behaviour of the
argument.

Given a sector $ S=\{w\colon 0<|w|<\infty, \ \theta_1<\arg
w<\theta_2\} $ of opening $ \alpha = \theta_2-\theta_1 $ and an
entire function $ f $ (everywhere below it will be assumed
non-constant), consider the relative area of the preimage
$ f^{-1}S $:
\[
A(r, S, f) = \frac{\A (f^{-1}S\cap r\D)}{\A (r\D)}
\]
where $r\D = \{z\colon |z|<r\}$.

\begin{thm}\label{thm1}
Let $f$ be an entire function of finite positive order $\rho$.
Then there exist arbitrarily large $r$ such that for every sector
$S$ of opening $\alpha$, we have
\[
A(r, S, f) \ge \frac{\alpha \cdot \kappa(\rho)}{2\pi}
\]
where $\kappa (\rho)>0$ depends only on $\rho$, and
$\kappa(\rho)\sim \text{const}\,\rho^{-1}$ for $\rho\to \infty$.
\end{thm}

Plausibly, the statement of Theorem~\ref{thm1} can be
complemented by  $\kappa (\rho) \to 1$ as $\rho\to 0$.
For entire functions of order zero the equidistribution pattern is
more visible since they behave like monomials $z^n$ on a sequence
of wide annuli.
\begin{thm}\label{thm2}
Let $f$ be an entire function of order zero. Then, given
$\varepsilon>0$, there exist arbitrary large values $r$ such that for
every sector $S$ of opening $\alpha$,
\begin{equation}\label{eq1.0}
\frac{\alpha-\varepsilon}{2\pi} \le A(r, S, f) \le
\frac{\alpha+\varepsilon}{2\pi}\,.
\end{equation}
\end{thm}

The upper bound in \eqref{eq1.0} follows from the lower bound
applied to the complementary sector $\C \setminus S$.
\begin{cor}
Suppose $f$ is an entire function of order zero. Then for any
sector $S$ of opening $\alpha$,
\[
\liminf_{r\to\infty} A(r, S, f) \le \frac{\alpha}{2\pi} \le
\limsup_{r\to\infty} A(r, S, f)\,.
\]
\end{cor}

Note that for functions of positive order the lower limit
\[\liminf_{r\to\infty} A(r, S, f)\] is not obliged to be small when
$\alpha$ is small. If the order $\rho>1/2$,
consider an entire function $f$ that tends to $1$ as $z\to\infty$
uniformly
within some angle, and take $S=\{w\colon -\alpha/2 < \text{arg}(w)
<\alpha/2 \}$.
If the order $\rho\le 1/2$, the following example was suggested by
Alexander Fryntov: take
\[
f(z) = 1 + \prod_{k=1}^\infty \left( 1 - \frac{z}{T^k}
\right)^{[T^{k\rho}]}\,,
\]
where $T$ is a sufficiently large positive number.

The next result gives a non-asymptotic version
of the coarse equidistribution principle. Define the
{\em doubling exponent} $ \beta(D,f) $ of an
analytic function $f$ on the disc $D$ as
\[
\beta (D, f) = \log \frac{\sup_D |f|}{\max_{\frac12 D}|f|}\,.
\]
It measures a certain complexity of the function $f$, cf \cite{NY, RY}.
For example,
if $f$ is a polynomial of degree $d$, then it is not difficult to
see that $ \beta (D, f) \le Cd$ for any disc $D$ in $\C$.
\begin{thm}\label{thm4}
Let $f$ be a non-zero analytic function in the unit disc $\D$,
$f(0)=0$. Then for any sector $S$ of opening $\alpha$ (with vertex
at the origin)
\[
A(r, S, f) \ge \frac{c\alpha}{\log\beta^*(\D, f)}\,,
\]
where $c$ is a positive numerical constant, and $\beta^* =
\max(\beta, 2)$.
\end{thm}

A special case of Theorem~\ref{thm4} with $S=\{w\colon \text{Re}
w>0 \}$ appeared in our recent work with L.~Polterovich
\cite[Theorem~2.2]{NPS}. It was preceded by a qualitative
compactness lemma proved by Nadirashvili in \cite{Nad}.

Concluding this introduction, we mention a curious resemblance
between Theorem~\ref{thm4} and a result of Marshall and
Smith~\cite{MS} that says that for any univalent analytic
function $f$ in $\D$ and any sector $S$ of opening $\alpha$
\begin{equation}\label{eq1.MS}
\iint_{f^{-1}S} |f|\, d\A \ge \kappa (\alpha) \iint_{\D} |f|\,
d\A\,,
\end{equation}
where $\kappa$ depends only on $\alpha$. It remains an open
question whether \eqref{eq1.MS} persists for arbitrary analytic
functions $f$ in $\D$ vanishing at the origin.

\medskip\par\noindent{\sc Organization of the paper.} In Section~2, we
introduce a characteristic $\Omega (r, f)$ which measures
oscillation of $\arg f$ on concentric circles and increases with
$r$. Then we formulate our Main Lemma and prove
Theorems~\ref{thm1} and~\ref{thm4}. In Section~3, we prove the
Main Lemma. In Section~4, we prove Theorem~\ref{thm2}. This part is
independent from the previous sections.

\medskip\par\noindent{\sc Convention.} Notation $ A \lesssim B$
means that $ A \le cB$ where $c$ is a positive numerical constant.

\medskip\par\noindent{\sc Acknowledgement. } This note grew out of our
joint work with Leonid Polterovich \cite{NPS}.
We thank him and Alexander Fryntov for numerous inspiring discussions.

\section{Oscillation of argument and the main lemma}

Suppose that the function $f$ is analytic in the disc $R\D$ and
does not vanish on the circle $r\T$, $0<r<R$. Consider all arcs
$L\subset r\T$ traveled counterclockwise (including the entire
circumference $r\T$ viewed as an arc whose end and beginning coincide).
Put
$$
\omega(r, f) :=\max_{L\subset r\T}\Delta_L\arg f
$$
where $\Delta_L\arg f$ is the increment of the argument of $f$
over $L$. The function
$r\mapsto \omega (r, f)$ is not necessarily monotone. To fix this
drawback, we slightly modify the definition and define a monotone
function $\Omega (r, f) $ that is close to $\omega (r, f)$.
By $ \Zf $ we denote the zero set of $ f $ (the zeroes are counted
with their multiplicities). Let $n(r, f) = \# \Zf \cap r\clD $,
$|\Zf| = \{|\zeta|\colon\zeta\in\Zf\} $.

Given $r\in (0, R)\setminus |\Zf|$, consider the factorization
\[
f(z) = e^{g_r(z)}\, \prod_{\zeta\in\Zf\cap r\D} (z-\zeta)
\]
and take
\[
\Omega (r, f) = 2\pi n(r, f) +
\underset{r\T}{\text{osc}}(\text{Im} g_r)\,,
\]
where
\[
\underset{r\T}{\text{osc}}(h) = \max_{r\T} h - \min_{r\T} h
\]
is the oscillation of the function $h$ on the circle $r\T$.
Below we present several properties of the characteristic $\Omega
(r, f)$.

\medskip\par\noindent 1. {\em The function $\Omega (r, f)$ increases
with $r$.} Indeed, if $r_1<r_2$ and there are no zeroes of $f$ in
the annulus $\{r_1\le |z| \le r_2\}$, then $\text{Im} g_{r_1}$ and
$\text{Im} g_{r_1}$ are traces of the same harmonic function on
different circles, and by the maximum principle, the oscillation of
any harmonic function on the circle increases when the radius
increases. Now, let us see what happens when $r$ runs through
$r_0$ such that $f$ vanishes on $r_0\T$. If $\zeta_0\in r_0\T$ is a
zero of $f$ of multiplicity $m$, then we need to add $m$ to $n(r,
f)$, and subtract $m\arg (z-\zeta_0)$ from $\text{Im} g_r$. Since
\[
\lim_{\varepsilon\to 0}\underset{_{(r_0-\varepsilon)\T}}{\text{osc}} \arg
(z-\zeta_0) =
\pi\,,
\]
we see that in this case
\[
\underset{(r_0+\varepsilon)\T}{\text{osc}}(\text{Im}\,g_{r_0+\varepsilon})
-
\underset{(r_0-\varepsilon)\T}{\text{osc}}(\text{Im}\,g_{r_0-\varepsilon})
\ge  -m \pi \,.
\]
Thus  $\Omega (r_0+0, f) > \Omega(r_0-0, f) $.

\medskip\par\noindent 2. {\em For all $r\in (0, R)\setminus |\Zf|$,
\[
\frac12 \Omega (r, f) \le \omega (r, f) \le \Omega (r, f)\,.
\]
}
The upper bound follows by definitions of $\omega$ and $\Omega$.
The lower bound is also easy: if $2\pi n(r, f) \ge \frac12 \Omega
(r, f)$, then
\[
\omega (r, f) \ge 2\pi n(r, f) \ge \frac12 \Omega (r, f)\,.
\]
If $2\pi n(r, f) < \frac12 \Omega (r, f)$, then the oscillation of
$\text{Im} g_r$ on $r\T$ is larger than $\frac12 \Omega (r, f)$.
Consider the arc on $r\T$ that runs counterclockwise from the
minimum to the maximum of $\text{Im} g_r$. The increment of $\arg
f $ on this arc cannot be smaller than the oscillation of
$\text{Im} g_r$; i.e., than $\frac12 \Omega (r, f)$.

\medskip\par\noindent 3. {\em Let
\[
\beta (r, f) = \log M(r, f) - \log M(\tfrac12 r, f)\,.
\]
Then
\[
\Omega (\tfrac12 r, f) \lesssim \beta^*(r, f) \qquad \text{and}
\qquad \beta (\tfrac12 r, f) \lesssim \Omega^*(r, f)\,,
\]
where $a^*:=\max(a,2)$.}

This inequalities go back to Gelfond~\cite{G} and
Hellerstein-Korevaar~\cite{HK}. We shall use only the first bound
whose proof can be found in \cite{NPS}. This bound immediately
yields the following property

\medskip\par\noindent 4. {\em Suppose $f$ is an entire function. Then
\[
\Omega (r, f) \lesssim \log M(2r, f)\,, \qquad r\ge r_0(f)\,.
\]
}

\medskip
For $L\subset [0, \infty)$, we set $K_L=\{z\colon |z|\in L\}$. The
following lemma plays a central role:
\begin{lemma}\label{lemma_main}
Suppose that the analytic function $f$ on $\clD$ and $t\in (0, 1)$
are such that
\[
\inf_{[t, 1]} \omega (r, f) \ge 2\pi\,,
\]
and
\[
\Omega (t, f) \ge \frac12 \Omega (1, f)\,.
\]
Then
\[
\A (f^{-1}S \cap K_{[t, 1]}) \gtrsim \alpha (1-t)^2\,.
\]
\end{lemma}

Now, we deduce Theorems~\ref{thm1} and~\ref{thm4} from this lemma.
The lemma will be proven in the next section.

\medskip\par\noindent{\em Proof of Theorem~\ref{thm1}:}
Since the function $f$ has order $\rho$, we can find arbitrarily
large $r$ such that
\[
\Omega (2r, f) \le 2^{2\rho} \Omega (r, f)
\]
and
\[
\omega (t, f) > 2\pi
\]
for $t\ge r$. Assume that $\rho\ge 2$ and split the interval $[r, 2r]$
into
$[5\rho]$ equal parts by points $r_0=r$, $r_1$, ..., $r_{[5\rho]}
= 2r$. Note that the inequality
\[
\Omega (r_{j+1}, f) \le 2^4 \Omega (r_j, f)
\]
holds for at least half of the indices $j=1$, ..., $[5\rho]$. For
these $j$, by Lemma~\ref{lemma_main}, the relative area of the set
$f^{-1}S\cap K_{[r_j, r_{j+1}]}$ with respect to $K_{[r_j, r_{j+1}]}$
is $\gtrsim \alpha \rho^{-1}$. Hence,
\[
A(r, S, f) \gtrsim \frac{\alpha}{\rho}\,,
\]
and we are done. \hfill $\Box$

\medskip\par\noindent{\em Proof of Theorem~\ref{thm4}:} Choose
$k\in\mathbb N$ such that
\[
2^k \pi \le \Omega (r, f) < 2^{k+1}\pi
\]
(recall that $f(0)=0$, thus $\Omega (1, f) \ge 2\pi$). Choose
$0=r_0<r_1<\,...\,<r_k\le 1$ so that $\Omega (r_j-0, f) \le 2^j\pi
\le \Omega (r_j+0, f)$, and set $r_{k+1}=1$.
Applying (properly scaled) Lemma~\ref{lemma_main} to the annuli
$K_j=K_{[r_j, r_{j+1}]}$, we get
\[
\A (f^{-1}\cap K_j) \gtrsim \alpha (r_{j+1} - r_j)^2\,.
\]
By Cauchy's inequality,
\[
\sum_{j=0}^k (r_{j+1}-r_j)^2 \ge \frac1{k+1}\,.
\]
Therefore,
\[
\A (f^{-1}S \cap \D) \ge \sum_{j=1}^k \A (f^{-1}S \cap K_j)
\gtrsim \frac{\alpha}{k}\,,
\]
completing the proof. \hfill $\Box$


\section{Proof of the main lemma}

Let, as above, $S=\{w\colon 0<|w|<\infty, \theta_1<\arg w
<\theta_2\}$, $\theta_2-\theta_1=\alpha$. Fix $r\in [t,1]\setminus
|\Zf|$, and introduce two types of `traversing arcs' on $r\T$:
$T$-arcs and $S$-arcs. An open arc $J\subset r\T$ is called a
$T$-arc, if a continuous branch of $\arg f$ maps $J$ onto an
interval $(\theta_1 + 2\pi m, \theta_1 + 2\pi (m+1))$ for some
$m\in \mathbb Z$. Each $T$-arc $J$ contains a traversing $S$-arc
$I$ which is mapped by the same branch of $\arg f$ onto $(\theta_1
+ 2\pi m, \theta_1+\alpha + 2\pi m)$. For each $r\in
[t,1]\setminus |\Zf|$, there are at least
\[
M = \left[ \frac1{2\pi} \inf_{[t,1]\setminus |\Zf|} \omega (r, f)
\right] \ge 1
\]
disjoint traversing $T$-arcs. We choose $M$ of them and discard
the rest.

Let $t_1 = \frac12 (1+t)$, $K=K_{[t, t_1]}$, and let $E$ be the
union of all $S$-arcs in $K$. We need to estimate from below the
area of $E$. Start with the argument used in \cite{NPS}. For each
$S$-arc $I\subset r\T$,
\[
\int_I |\nabla \arg f|\, |dz| \ge \alpha\,.
\]
Therefore,
\[
\int_{E\cap r\T} |\nabla \arg f|\, |dz| \ge \alpha M\,.
\]
Integrating by $r\in [t, t_1]$, we get
\begin{equation}\label{eq3.1}
\iint_E |\nabla \arg f|\, d\A \ge \frac{\alpha (1-t)M}2\,.
\end{equation}
Now we shall try to estimate from above the double integral on
the left-hand side.

Factoring
\[
f(z) = e^{g(z)} \prod_{\zeta\in\Zf} (z-\zeta)\,,
\]
we get
\[
|\nabla \arg f| \le |\nabla \text{Im}\, g| + \sum_{\zeta\in\Zf}
\frac1{|z-\zeta|}\,,
\]
and
\[
\iint_E |\nabla \arg f|\, d\A  \le \max_{t_1\clD}
|\nabla\text{Im}\, g| \cdot \A(E) + \# \Zf \cdot
\sup_{\zeta\in\Zf} \iint_E \frac{d\A(z)}{|z-\zeta|}\,.
\]
Estimate the terms on the right-hand side. We have $\Omega (1, f)
\le 16\pi M $,
\[
\# \Zf \le \frac1{2\pi} \Omega (1, f) \le 8M\,,
\]
and
\[
\underset{\T}{\text{osc}} (\text{Im}\, g) \le \Omega (1, f) \le
16\pi M\,.
\]
Since the function $\text{Im}\, g$ is harmonic in $\D$, we obtain
\[
\max_{t_1\clD} |\nabla \text{Im}\, g| \le \frac1{1-t_1}\cdot
\underset{\T}{\text{osc}} (\text{Im}\, g) \le \frac{32\pi
M}{1-t}\,,
\]
and finally
\begin{equation}\label{eq3.2}
\iint_E |\nabla \arg f|\, d\A  \lesssim M \left( \frac{\A(E)}{1-t}
+ \sup_{\zeta\in\Zf} \iint_E \frac{d\A(z)}{|z-\zeta|} \right) \,.
\end{equation}
It remains to estimate the double integral on the right-hand side.
Till that point we followed the strategy from \cite{NPS}. A
straightforward bound
\begin{equation}\label{eq3.3}
\iint_E \frac{d\A(z)}{|z-\zeta|} \le 2\sqrt{\pi \A(E)}
\end{equation}
used in there is not sufficient anymore:\footnote{However, it will be
employed below during an auxiliary step.} it leads only to the
estimate
\[
\A (f^{-1}S \cap K) \gtrsim \alpha^2 (1-t)^2\,.
\]
We try to get something better taking into account the structure
of the set $E$ (recall that $E\cap r\T$ is always a union of $M$
disjoint $S$-arcs). For this purpose, we reduce the general case
to the one when all $S$-arcs are short, the $T$-arcs containing
them are not very short (that is, the $S$-arcs are `well
separated'), and the zero set $\Zf$ is not too close to $E$.

First, we sort the $S$-arcs. We call an $S$-arc $I$ a {\em short}
one, if
\begin{equation}\label{eq3.short}
|I| \le \frac{\alpha \eta (1-t)}{M}\,,
\end{equation}
where a small positive numerical constant $\eta$ will be chosen later.
Otherwise, we say that $I$ is not short. By $M_s(r)$ we denote the
number of short $S$-arcs on $r\T$. Let $E_{n.s.}$ be the union of
all non-short arcs in $K$. Clearly,
\begin{equation}\label{eq3.4}
\A (E_{n.s.}) \ge \frac{\alpha \eta (1-t)}{M}\, \int_t^{t_1} (M -
M_s(r) )\, dr\,.
\end{equation}

Now consider the short $S$-arcs in $K$. In fact, we do not need all
of them. Let $E_s^*$ be the union of
all $S$-arcs $I$ in $K$ satisfying the following three conditions:

\smallskip\par\noindent (a) $I$ is short (i.e. \eqref{eq3.short} holds);

\smallskip\par\noindent (b) the corresponding $T$-arc $J\supset I$
is not very short:
\[
|J| \ge \frac{\delta (1-t)}{M}\,,
\]
where a small positive numerical constant $\delta$ will be chosen later;

\smallskip\par\noindent (c) if $I\subset r\T$, then
\[
\text{dist}(r, |\Zf|) \ge \frac{\delta (1-t)}{M}\,.
\]

We will show that under appropriate choice of small parameters $\delta$
and
$\eta$,
\begin{equation}\label{eq.*}
\A (E_s^*) \gtrsim \frac{\alpha (1-t)}{M} \left( \int_t^{t_1} M_s(r)\, dr
-
\frac{(1-t)M}5 \right)\,.
\end{equation}
Then recalling \eqref{eq3.4}, we get the assertion of the main lemma.
If $ \int_t^{t_1} M_s(r) \, dr \ge \frac{(1-t)M}{5} $, then
\begin{eqnarray*}
\A (f^{-1}S \cap K) &\ge& \A (E_{n.s.}) + \A (E_s^*) \\
&\gtrsim& \frac{\alpha (1-t)}{M} \left( (t_1-t) M -
\frac{(1-t)M}{5} \right) \gtrsim \alpha (1-t)^2\,.
\end{eqnarray*}
If $\int_t^{t_1} M_s(r) \, dr < \frac{(1-t)M}{5}$,
then we simply discard the short $S$-arcs:
\begin{eqnarray*}
\A (f^{-1}S \cap K) &\ge& \A (E_{n.s.}) \\
&\stackrel{\eqref{eq3.4}}\ge& \frac{\alpha \eta (1-t)}{M}
\left( (t_1 - t)M - \frac{(1-t)M}{5} \right)
\gtrsim \alpha (1-t)^2\,.
\end{eqnarray*}

Now we start proving \eqref{eq.*}.
Let
\[
\mathcal E = \left\{r\in [t, t_1]\colon \text{dist}(r, |\Zf|) <
\tfrac{\delta (1-t)}{M} \right\}
\]
be an exceptional set of radii, and let $m(r)$ be the number of
{\em very short} $T$-arcs $J$ such that
\[
|J| \le \frac{\delta (1-t)}{M}\,.
\]
Then, as above, for $r\notin\mathcal E$,
\[
\int_{E_s^* \cap r\T} |\nabla \arg f|\, |dz| \ge \alpha (M_s (r) -
m(r) )\,,
\]
and
\begin{eqnarray*}
\iint_{E_s^*} |\nabla \arg f|\, d\A &\ge& \alpha \int_{[t,
t_1]\setminus \mathcal E} (M_s(r)-m(r))\, dr \\
&\ge& \alpha \left( \int_t^{t_1} M_s(r)\, dr - \int_t^{t_1} m(r)\,
dr - M |\mathcal E| \right)\,.
\end{eqnarray*}
The next two claims show that the second and the third terms on
the right-hand side are relatively small, provided that $\delta$ is
sufficiently small.

\begin{claim}\label{claimA} We have
\[
\int_t^{t_1} m(r)\, dr \lesssim \delta(1-t)M\,.
\]
\end{claim}

\smallskip\par\noindent{\em Proof:} Let $G$ be the union of all very
short $T$-arcs in $K$. Then, as above,
\[
\int_{G\cap r\T} |\nabla \arg f|\, |dz| \ge 2\pi m(r)\,,
\]
and
\[
\iint_G |\nabla \arg f|\, d\A \ge 2\pi \int_t^{t_1} m(r)\, dr\,.
\]
On the other hand, a counterpart of \eqref{eq3.2} together with
estimate \eqref{eq3.3} give us
\[
\iint_G |\nabla \arg f|\, d\A \lesssim M \left( \frac{\A (G)
}{1-t} + \sqrt{\A (G)}\right)\,.
\]
If the second term on the right hand side is larger than the first one,
then
we get
\begin{equation}\label{eq.new}
\A (G) \gtrsim \left( \frac1{M} \int_t^{t_1} m(r)\, dr\right)^2\,.
\end{equation}
If the first term is larger, then $\A (G) \ge (1-t)^2$ and we again arrive
at
\eqref{eq.new}

Since $G$ consists of very short $T$-arcs, we have
\[
\A (G) \le \frac{\delta (1-t)}{M}\, \int_t^{t_1} m(r)\, dr\,.
\]
Hence,
\[
\frac1{M} \int_t^{t_1} m(r)\, dr \lesssim \delta (1-t)\,,
\]
proving the claim. \hfill $\Box$

\begin{claim}\label{claimB} We have
\[
|\mathcal E| \lesssim \delta (1-t)\,.
\]
\end{claim}

\smallskip\par\noindent{\em Proof:} Since $\# \Zf \lesssim M$, this
follows from definition of $\mathcal E$. \hfill $\Box$

\medskip Using these claims, we choose $\delta$ so small that
\[
\int_t^{t_1} m(r)\, dr + M|\mathcal E| \le \frac{(1-t)M}{10}\,.
\]
Then
\begin{equation}\label{eq3.5}
\iint_{E_s^*} |\nabla \arg f|\, d\A \ge \alpha \left( \int_t^{t_1}
M_s(r)\, dr - \frac{(1-t)M}{10} \right)\,.
\end{equation}
We are ready to make the final step: to estimate from above the
integral on the left-hand side of \eqref{eq3.5}. As above (cf.
\eqref{eq3.2}\,),
\begin{equation}\label{eq.^}
\iint_{E_s^*} |\nabla \arg f|\, d\A  \lesssim M \left(
\frac{\A(E_s^*)}{1-t} + \sup_{\zeta\in\Zf} \iint_{E_s^*}
\frac{d\A(z)}{|z-\zeta|} \right) \,.
\end{equation}
The next claim bounds the double integral on the right-hand side:
\begin{claim}\label{claimC}
Let $F\subset K_{[t, 1]} $ be a closed set, $t\ge \frac12$.
Suppose that there exists
$s\in (0, \frac{1-t}2)$ such that, for each $r\in (t, 1)$,
the set
$F(r)=F\cap r\T$ is a union of disjoint arcs $I$ of length
\[
|I|\le \beta s.
\]
Further, assume that each arc $I$ is contained in a bigger arc $J$,
the arcs $J$ are pairwise disjoint,
\[
|J| \ge s\,,
\]
and the total number of arcs is $\lesssim \frac{1-t}{s}$.
Given $\xi\in [t, 1]$, denote
\[
F_\xi = F\setminus K_{[\xi-s, \xi+s]}\,.
\]
Then
\[
\iint_{F_\xi} \frac{d\A(z)}{|z-\xi|} \lesssim \beta(1-t)\,.
\]
\end{claim}

Note that under the assumptions of this claim, $\A (F_\xi) \lesssim
\beta (1-t)^2$, and estimate \eqref{eq3.3} gives us only
\[
\iint_{F_\xi} \frac{d \A(z)}{|z-\xi|} \lesssim \sqrt{\beta}
\cdot(1-t)\,.
\]

\par\noindent{\em Proof of Claim~\ref{claimC}:} Fix $r\in (t, 1)$,
$|r-\xi|\ge s$, and consider the integral
\[
\int_{F(r)} \frac{|dz|}{|z-\xi|}\,.
\]
First, estimate the contribution of the components $I$ of $F(r)$ that
intersect the arc $\{re^{i\theta}\colon |\theta|\le |r-\xi|\}$. For each
arc $I$,
\[
\int_I \frac{|dz|}{|z-\xi|} \le
\frac{|I|}{|r-\xi|} \le \frac{\beta s }{|r-\xi|}\,.
\]
The number of such arcs $I$ is $\lesssim \frac{|r-\xi|}{s}$.
Therefore, the total contribution of these arcs is $\lesssim
\beta (1-t)$.

Now consider the arcs $I$ that do not intersect the arc
$\{re^{i\theta}\colon |\theta|\le |r-\xi|\}$. It suffices
to consider only the arcs $I$ lying in the upper semi-circle. For
these arcs, $|re^{i\theta}-\xi| \gtrsim  \theta $.
We enumerate the arcs $I$ counterclockwise by index $j$,
$1\le j \lesssim \frac{1-t}{s}$.
Then the contribution of the $j$-th arc is
\[
\int_{I_j} \frac{|dz|}{|z-\xi|} \lesssim
\int_{|r-\xi| + (j-1)s}^{|r-\xi| + (j-1)s +\beta s} \frac{d\theta}{\theta}
\le
\frac{\beta}{ \frac{|r-\xi|}{s} + j-1}\,.
\]
Summing over $j$, we get the bound
\[
\beta \left( \log \frac{1-t}{|r-\xi|} + \text{Const}\right)\,.
\]
Integrating this bound by $r$ from $t$ to $1$, we see that the
contribution of these arcs is $\lesssim \beta (1-t)$ as well. \hfill
$\Box$

\medskip At last, we are able to get estimate \eqref{eq.*} for
$ \A(E_s^*) $, and thus to finish the proof of Lemma~\ref{lemma_main}.
Without loss of generality, we assume that $t\ge \frac12$.
We apply the claim to the integral on
the right hand side of \eqref{eq.^} with
\[
\beta = \frac{\alpha \eta}{\delta}\,, \qquad s =
\frac{\delta (1-t)}{M}
\]
(recall that the parameter $\delta$ already has been fixed,
but $\eta$ has not been chosen yet). We obtain
\begin{equation}\label{eq3.6}
\iint_{E_s^*} |\nabla \arg f|\, d\A  \lesssim M \left(
\frac{\A(E_s^*)}{1-t} + \frac{\alpha \eta (1-t)}{\delta} \right)
\,.
\end{equation}
Juxtaposing estimates \eqref{eq3.5} and \eqref{eq3.6}, we get
\[
\A (E_s^*) \gtrsim \frac{\alpha (1-t)}{M} \left( \int_t^{t_1}
M_s(r)\, dr - \frac{(1-t)M}{10} \right) - \frac{\alpha \eta
(1-t)^2}{\delta}\,.
\]
It remains to choose $\eta$ so small that the right-hand side is
\[
\gtrsim \frac{\alpha (1-t)}{M} \left( \int_t^{t_1} M_s(r)\, dr -
\frac{(1-t)M}5 \right)\,.
\]
This gives us \eqref{eq.*} and completes the proof of the main
lemma. \hfill $\Box$

\section{Functions of order zero}

Here we prove Theorem~\ref{thm2}. Without loss of generality,
$f(0)\ne 0$ and $f$ is not a polynomial. Then, up to a constant
factor that is irrelevant here, we have
$$
f(z)=\prod_{\zeta\in\mathcal Z_f}\left(1-\frac{z}{\zeta}\right)\,.
$$
Fix an $\varepsilon>0$. Then take a very small $\delta >0$ to be
chosen later. Choose $r_\delta$ to be the radius at which the
ratio
$$\frac{n(r, f)}{r^\delta}$$ attains its maximum. Note that
$r_\delta\to+\infty$ as $\delta \to 0+$. Let $M=n(r_\delta, f)$.
Let $U$ be a huge constant (depending only on $\varepsilon$) to be
chosen later. We claim that the disk $R\D $ with $R=U^2r_\delta$
satisfies the equidistribution property of the theorem if $\delta
$ is small enough.

Indeed, consider the annulus $K:=\{z\colon Ur_\delta<|z|<R\}$. For
every $r\in(Ur_\delta,R)$ the set $f^{-1}S \cap r\T$ contains at
least $M$ disjoint traversing $S$-arcs and, thereby,
$$
\int_{f^{-1}S \cap r\T}|\nabla \arg f(z)|\,|dz|\ge M\alpha\,.
$$
Therefore,
\begin{equation}\label{eq4.1}
\iint_{E}|z|\,|\nabla \arg f(z)|\,d\A(z)\ge M\alpha \int_{R/U}^R
r\,dr= M\alpha \frac{R^2}{2}(1-U^{-2})\,,
\end{equation}
where $E=f^{-1}S \cap K$.

Now, we estimate the double integral on the left hand side from
above. Write $f=f_1 \cdot f_2 \cdot f_3$, where
\begin{align*}
f_1 &= \prod_{\zeta\in \Zf \,,\,|\zeta|\le r_\delta}
\left(1-\frac{z}{\zeta}\right)
\\
f_2&=\prod_{\zeta\in \Zf\,,\,r_\delta<|\zeta|\le U^3r_\delta}
\left(1-\frac{z}{\zeta}\right)
\\
f_3&=\prod_{\zeta\in \Zf\,,\,|\zeta|>U^3 r_\delta}
\left(1-\frac{z}{\zeta}\right)\,.
\end{align*}
For $z\in K$ we have
$$
|\nabla \arg f_1(z)|\le \sum_{\zeta\in \mathcal Z_f\,,\,|\zeta|\le
r_\delta} \frac{1}{|z-\zeta|}\le \frac{U}{U-1}\frac{M}{|z|}.
$$
Also
\begin{eqnarray*}
|\nabla \arg f_3(z)| &\le& \sum_{\zeta\in \mathcal
Z_f\,,\,|\zeta|>U^3r_\delta} \frac{1}{|z-\zeta|} \le \frac{U}{U-1}
\sum_{\zeta\in \mathcal Z_f\,,\,|\zeta|>
U^3r_\delta}\frac{1}{|\zeta|}
\\
&=& \frac{U}{U-1}\sum_{j\ge 3} \ \sum_{\zeta\in
\Zf,\,U^jr_\delta<|\zeta| \le U^{j+1}r_\delta}\frac{1}{|\zeta|}
\\
&\le& \frac{U}{U-1}\sum_{j\ge 3} \frac{1}{U^j
r_\delta}(U^{(j+1)\delta}-1)M
\\
&=& \frac{U}{U-1}\frac{1}{U^2r_\delta} \sum_{j\ge
1}\frac{1}{U^j}(U^{(j+3)\delta}-1)M \le
\frac{U}{U-1}\sigma(U,\delta)\frac{M}{|z|}
\end{eqnarray*}
where
$$
\sigma(U,\delta)=\sum_{j\ge 1}\frac{1}{U^j}(U^{(j+3)\delta}-1) \to
0 \qquad \text{as } \delta\to 0 \quad\text{for any fixed }U>1\,.
$$
Therefore,
\[
|\nabla \arg f_1(z)|+ |\nabla \arg f_3(z)|\le
(1+\gamma(U,\delta))\frac{M}{|z|}\,,
\]
where $\gamma(U,\delta)$ can be made arbitrarily small if $U$ is
large enough and $\delta$ is small enough. Note that the number of
zeroes in $f_2$ does not exceed $(U^{3\delta}-1)M$. Hence,
\begin{eqnarray*}
\iint_E |z|\, |\nabla \arg f_2(z)|\, d\A(z) &\le&
(U^{3\delta}-1)MR \iint_{|z|\le R} \frac{d\A(z)}{|z-\zeta|} \\
&\le& (U^{3\delta}-1)M \cdot 2\pi R^2\,.
\end{eqnarray*}
Thus
\begin{equation}\label{eq4.3}
\iint_E |z|\,|\nabla \arg f(z)|\,d\A(z)\le
(1+\gamma(U,\delta))M\A(E) + (U^{3\delta}-1)M\cdot 2\pi R^2\,.
\end{equation}

The rest is clear. We choose $U$ so large such that
$U^{-2}<\frac{\varepsilon}4$. Then we choose $\delta$ so small that
\[
\gamma (U, \delta) < \frac{\varepsilon}4\,, \qquad \text{and} \qquad
U^{3\delta} -1 < \frac{\varepsilon}4\,.
\]
Juxtaposing \eqref{eq4.1} and \eqref{eq4.3}, cancelling $M$, and
taking into account the choice of $U$ and $\delta$, we get the
result. \hfill $\Box$

\bigskip

\begin{tabbing}
F\"edor Nazarov \qquad\qquad\qquad\qquad\qquad\qquad\qquad\qquad\=
Mikhail Sodin\\
Department of Mathematics \> School of Mathematics\\
Michigan State University \> Tel Aviv University\\
East Lansing, MI 48824 \> Tel Aviv 69978\\
USA \> Israel\\
\texttt{\small fedja@math.msu.edu} \>
\texttt{\small sodin@post.tau.ac.il}
\end{tabbing}


\begin{thebibliography}{}


\bibitem{G} {\sc A.~Gelfond}, \"Uber die harmonischen Funktionen,
Trav. Inst. Stekloff {\bf 5} (1934), 149-158.

\bibitem{HK} {\sc S. Hellerstein and J. Korevaar},
The real values of an entire function. Bull. Amer. Math. Soc. {\bf
70} (1964), 608--610.

\bibitem{MS} {\sc D. Marshall and W. Smith},
The angular distribution of mass by Bergman functions, Rev. Mat.
Iberoamericana {\bf 15} (1999), 93--116.


\bibitem{Nad} {\sc N. Nadirashvili},
Metric properties of eigenfunctions of the Laplace operator on
manifolds, Ann. Inst. Fourier {\bf 41} (1991), 259--265.

\bibitem{NPS} {\sc F. Nazarov, L. Polterovich, and M. Sodin},
Sign and area in nodal geometry of Laplace eigenfunctions,
\texttt{arXiv math.AP/0402412}

\bibitem{NY} {\sc D. Novikov and S. Yakovenko},
A complex analogue of the Rolle theorem and polynomial envelopes of
irreducible differential equations in the complex domain,
J. London Math. Soc. (2) {\bf 56} (1997), 305--319.


\bibitem{RY} {\sc N. Roytwarf and Y. Yomdin},
Bernstein classes, Ann. Inst. Fourier {\bf 47} (1997), 825--858.



\end{thebibliography}
\end{document}